\documentclass{amsart}
\usepackage{amsmath,amscd,amsthm,amsfonts}

\theoremstyle{plain}
\newtheorem{theorem}                 {Theorem}      [section]
\newtheorem{proposition}  [theorem]  {Proposition}

\newtheorem{lemma}        [theorem]  {Lemma}
\newtheorem{conjecture}        [theorem]  {Conjecture}

\theoremstyle{definition}

\newtheorem{definition}   [theorem]  {Definition}

\numberwithin{equation}{section}

\begin{document}
\baselineskip 18pt \larger

\def \theo-intro#1#2 {\vskip .25cm\noindent{\bf Theorem #1\ }{\it #2}}

\newcommand{\trace}{\operatorname{trace}}

\def \dn{\mathbb D}
\def \nn{\mathbb N}
\def \zn{\mathbb Z}
\def \qn{\mathbb Q}
\def \rn{\mathbb R}
\def \cn{\mathbb C}
\def \hn{\mathbb H}
\def \P{\mathbb P}
\def \can{Ca}

\def \S{\mathcal S}
\def \A{\mathcal A}
\def \B{\mathcal B}
\def \C{\mathcal C}
\def \G{\mathcal G}
\def \F{\mathcal F}
\def \H{\mathcal H}
\def \I{\mathcal I}
\def \L{\mathcal L}
\def \M{\mathcal M}
\def \N{\mathcal N}
\def \Ol{\mathcal O}

\def \R{\mathcal R}
\def \V{\mathcal V}
\def \W{\mathcal W}

\def\Re{\mathfrak R\mathfrak e}
\def\Im{\mathfrak I\mathfrak m}
\def\Co{\mathfrak C\mathfrak o}
\def\Or{\mathfrak O\mathfrak r}

\def \ip #1#2{\langle #1,#2 \rangle}
\def \spl#1#2{( #1,#2 )}

\def \lb#1#2{[#1,#2]}

\def \pror#1{\rn P^{#1}}
\def \proc#1{\cn P^{#1}}
\def \proh#1{\hn P^{#1}}

\def \gras#1#2{G_{#1}(\cn^{#2})}

\def \b{\mathfrak{b}}
\def \g{\mathfrak{g}}
\def \h{\mathfrak{h}}
\def \k{\mathfrak{k}}
\def \m{\mathfrak{m}}
\def \p{\mathfrak{p}}
\def \q{\mathfrak{q}}
\def \r{\mathfrak{r}}
\def \un{\mathfrak{u}}

\def \GLR#1{\text{\bf GL}_{#1}(\rn)}
\def \GLRP#1{\text{\bf GL}^+_{#1}(\rn)}
\def \glr#1{\mathfrak{gl}_{#1}(\rn)}
\def \GLC#1{\text{\bf GL}_{#1}(\cn)}
\def \glc#1{\mathfrak{gl}_{#1}(\cn)}
\def \GLH#1{\text{\bf GL}_{#1}(\hn)}
\def \glh#1{\mathfrak{gl}_{#1}(\hn)}
\def \GLD#1{\text{\bf GL}_{#1}(\dn)}
\def \gld#1{\mathfrak{gl}_{#1}(\dn)}

\def \SLR#1{\text{\bf SL}_{#1}(\rn)}
\def \slr#1{\mathfrak{sl}_{#1}(\rn)}
\def \SLC#1{\text{\bf SL}_{#1}(\cn)}
\def \slc#1{\mathfrak{sl}_{#1}(\cn)}

\def \O#1{\text{\bf O}(#1)}
\def \SO#1{\text{\bf SO}(#1)}
\def \so#1{\mathfrak{so}(#1)}
\def \SOs#1{\text{\bf SO}^*(#1)}
\def \sos#1{\mathfrak{so}^*(#1)}
\def \SOO#1#2{\text{\bf SO}(#1,#2)}
\def \SOO0#1#2{\text{\bf SO}_0(#1,#2)}
\def \soo#1#2{\mathfrak{so}(#1,#2)}
\def \SOC#1{\text{\bf SO}(#1,\cn)}
\def \SOc#1{\text{\bf SO}(#1,\cn)}
\def \soc#1{\mathfrak{so}(#1,\cn)}

\def \SUS#1{\text{\bf SU}^*(#1)}
\def \sus#1{\mathfrak{su}^*(#1)}

\def \U#1{\text{\bf U}(#1)}
\def \u#1{\mathfrak{u}(#1)}
\def \US#1{\text{\bf U}^*(#1)}
\def \us#1{\mathfrak{u}^*(#1)}
\def \UU#1#2{\text{\bf U}(#1,#2)}
\def \uu#1#2{\mathfrak{u}(#1,#2)}
\def \SU#1{\text{\bf SU}(#1)}
\def \su#1{\mathfrak{su}(#1)}
\def \SUU#1#2{\text{\bf SU}(#1,#2)}
\def \suu#1#2{\mathfrak{su}(#1,#2)}

\def \Sp#1{\text{\bf Sp}(#1)}
\def \sp#1{\mathfrak{sp}(#1)}
\def \Spp#1#2{\text{\bf Sp}(#1,#2)}
\def \spp#1#2{\mathfrak{sp}(#1,#2)}
\def \SpR#1{\text{\bf Sp}(#1,\rn)}
\def \spR#1{\mathfrak{sp}(#1,\rn)}
\def \SpC#1{\text{\bf Sp}(#1,\cn)}
\def \spc#1{\mathfrak{sp}(#1,\cn)}

\def \d#1{\mathfrak{d}(#1)}
\def \s#1{\mathfrak{s}(#1)}
\def \sym#1{\text{Sym}(\rn^{#1})}
\def \symc#1{\text{Sym}(\cn^{#1})}

\def \gradh#1{\text{grad}_{\H}(#1 )}
\def \gradv#1{\text{grad}_{\V}(#1 )}

\def \nab#1#2{\hbox{$\nabla$\kern -.3em\lower 1.0 ex
    \hbox{$#1$}\kern -.1 em {$#2$}}}

\allowdisplaybreaks

\title{On the existence of harmonic morphisms\\
from certain symmetric spaces}

\author{Sigmundur Gudmundsson}
\author{Martin Svensson}

\thanks{The second author was supported by the Swedish Research Council (623-2004-2262)}



\address
{Mathematics, Faculty of Science, Lund University, Box 118, S-221
00 Lund, Sweden} \email{Sigmundur.Gudmundsson@math.lu.se}
\address
{Department of Pure Mathematics, University of Leeds, Leeds LS2
9JT, England} \email{M.Svensson@leeds.ac.uk}

\begin{abstract}
In this paper we give a positive answer to the open existence
problem for complex-valued harmonic morphisms from the non-compact
irreducible Riemannian symmetric spaces $\SLR n/\SO n$, $\SUS
{2n}/\Sp n$ and their compact duals $\SU n/\SO n$ and $\SU{2n}/\Sp
n$. Furthermore we prove the existence of globally defined,
complex-valued harmonic morphisms from any Riemannian symmetric
space of type IV.
\end{abstract}

\maketitle

\section{Introduction}

Harmonic morphisms are maps between Riemannian or semi-Riemannian
manifolds which pull back local harmonic functions on the codomain
to local harmonic functions on the domain. Equivalently, they may
be characterized as {\it harmonic maps} which satisfy the
additional condition of {\it horizontal (weak) conformality}.
Together, these two conditions form an over-determined, non-linear
system of partial differential equations, making the question of
the existence of harmonic morphisms interesting but very hard to
answer in general. Indeed, most metrics on a 3-dimensional domain
do not admit any non-constant solutions with values in a surface,
see \cite{Bai-Woo-1}.

In this paper we are mainly interested in maps with values in a
surface. In this case, the condition for a horizontally (weakly)
conformal map to be harmonic is equivalent to that of the map
having {\it minimal} regular fibres. Hence harmonic morphisms to
surfaces are useful tools to construct {\it minimal submanifolds}.
The equations for a map to a surface to be a harmonic morphism are
furthermore invariant under conformal changes of the metric on the
surface. Thus, at least for local studies, one can without loss of
generality assume that the codomain is the complex plane with its
standard metric.

It is known that in several cases, when the domain $(M,g)$ is an
irreducible Riemannian symmetric space, complex-valued solutions
to the problem do exist, see for example \cite{Gud-1}, \cite{Sve},
\cite{Gud-Sve-1} and \cite{Gud-Sve-2}. This has led the authors to
the following conjecture.

\begin{conjecture}\label{conj:existence}
Let $(M^m,g)$ be an irreducible Riemannian symmetric space of
dimension $m\ge 2$. For each point $p\in M$ there exists a
complex-valued harmonic morphism $\phi:U\to\cn$ defined on an open
neighbourhood $U$ of $p$. If the space $(M,g)$ is of non-compact
type then the domain $U$ can be chosen to be the whole of $M$.
\end{conjecture}

In this paper we introduce a {\bf new approach} to the problem and
employ this to prove the above conjecture in the cases when
$(M,g)$ is one of the non-compact irreducible Riemannian symmetric
spaces
$$\SLR n/\SO n, \qquad\qquad \SUS {2n}/\Sp n,$$ or their compact
dual spaces $$\SU n/\SO n, \qquad\qquad \SU{2n}/\Sp n.$$

In an earlier paper \cite{Gud-Sve-2} we constructed globally
defined complex-valued harmonic morphisms from those Riemannian
symmetric space $G^\cn/G$ of type IV, where $G$ is a simple
compact Lie group {\it admitting a Hermitian symmetric quotient}.
This proved Conjecture \ref{conj:existence} for the spaces
$$\SOC n/\SO n,\qquad \SLC n/\SU n,\qquad  \SpC n/\Sp n,$$
$$E_6^\cn/E_6,\qquad  E_7^\cn/E_7.$$
In the current paper we improve this result by showing that the
assumption that $G$ admits a Hermitian symmetric quotient is {\it
superfluous}. By this we prove Conjecture \ref{conj:existence} for
any Riemannian symmetric space of type IV and thereby add
$$E_8^\cn/E_8,\qquad F_4^\cn/F_4,\qquad G_2^\cn/G_2$$
to the list of spaces for which the statement is true.

We tacitly assume that all manifolds are connected
and that all objects such as manifolds, maps etc. are smooth, i.e.
in the $C^{\infty}$-category. For our notation concerning Lie groups
we refer to the comprehensive book \cite{Kna}.

\section{Harmonic Morphisms}

We are mainly interested in complex-valued harmonic morphisms
from Riemannian manifolds, but our methods involve maps
from the more general semi-Riemannian manifolds, see \cite{O-N}.

Let $M$ and $N$ be two manifolds of dimensions $m$ and $n$,
respectively. Then a semi-Riemannian metric $g$ on $M$ gives
rise to the notion of a
Laplacian on $(M,g)$ and real-valued harmonic functions
$f:(M,g)\to\rn$. This can be generalized to the concept of a
\emph{harmonic map} $\phi:(M,g)\to (N,h)$ between semi-Riemannian
manifolds, see \cite{Bai-Woo-book}.

\begin{definition}
A map $\phi:(M,g)\to (N,h)$ between semi-Riemannian manifolds is
called a {\it harmonic morphism} if, for any harmonic function
$f:U\to\rn$ defined on an open subset $U$ of $N$ with
$\phi^{-1}(U)$ non-empty, the composition
$f\circ\phi:\phi^{-1}(U)\to\rn$ is a harmonic function.
\end{definition}

The following characterization of harmonic morphisms between
semi-Riemannian manifolds is due to Fuglede, and generalizes the
corresponding well-known result of \cite{Fug-1,Ish} in the
Riemannian case.  See \cite{Bai-Woo-book} for the
definition of horizontal (weak) conformality.

\begin{theorem}\cite{Fug-2}
  A map $\phi:(M,g)\to (N,h)$ between semi-Rie\-mannian manifolds is a
  harmonic morphism if and only if it is both a harmonic map and
  horizontally (weakly) conformal.
\end{theorem}

The following result generalizes the corresponding well-known theorem
of Baird and Eells in the Riemannian case, see \cite{Bai-Eel}. It
gives the theory of harmonic morphisms a strong geometric flavour
and shows that the case when the codomain is a surface is
particularly interesting. Due to this result the conditions characterizing
harmonic morphisms are independent of conformal changes of the metric on the
surface.  For the definition of horizontal homothety we refer to
\cite{Bai-Woo-book}.

\begin{theorem}\cite{Gud-1}\label{theo:semi-B-E}
Let $\phi:(M^m,g)\to (N^n,h)$ be a horizontally conformal
submersion from a semi-Riemannian manifold $(M^m,g)$ to a
Riemannian manifold $(N^n,h)$. If
\begin{enumerate}
\item[(i)] $n=2$, then $\phi$ is harmonic if and only if $\phi$ has
minimal fibres,
\item[(ii)] $n\ge 3$, then two of the following conditions imply the other:
\begin{enumerate}
\item $\phi$ is a harmonic map,
\item $\phi$ has minimal fibres,
\item $\phi$ is horizontally homothetic.
\end{enumerate}
\end{enumerate}
\end{theorem}

\begin{proposition}\label{prop-lift}
Let $(\hat M,\hat g)$ be a semi-Riemannian manifold, $(M,g)$,
$(N,h)$ be Riemannian manifolds and $\pi:(\hat M,\hat g)\to(M,g)$
be a submersive harmonic morphism.  Furthermore let $\phi:(M,g)\to
(N,h)$ be a map and $\hat\phi:(\hat M,\hat g)\to(N,h)$ be the
composition $\hat\phi=\phi\circ\pi$.  Then $\phi$ is a harmonic
morphism if and only if $\hat \phi$ is a harmonic morphism.
\end{proposition}

\begin{proof}
Let $\lambda:\hat M\to\rn^+$ denote the dilation of the horizontally
conformal map $\pi:(\hat M,\hat g)\to (M,g)$.  If $f:U\to\rn$
is a function defined locally on $N$ then the composition law
for the tension field gives
\begin{eqnarray*}
\tau(f\circ\hat \phi)&=&\text{trace}\nabla
d(f\circ\phi)(d\pi,d\pi)+ d(f\circ\phi)(\tau(\pi))\\
&=&\lambda^2\tau(f\circ\phi)\circ\pi+d(f\circ\phi)(\tau(\pi))\\
&=&\lambda^2\tau(f\circ\phi)\circ\pi,
\end{eqnarray*}
because $\pi$ is horizontally conformal and harmonic.  The
statement then follows from the assumption that $\lambda^2>0$.
\end{proof}

In what follows we are mainly interested in complex-valued
functions $$\phi,\psi:(M,g)\to\cn$$ from semi-Riemannian manifolds.
In this situation the metric $g$ induces the complex-valued
Laplacian $\tau(\phi)$ and the gradient $\text{grad}(\phi)$ with
values in the complexified tangent bundle $T^{\cn}M$ of $M$.  We
extend the metric $g$ to be complex bilinear on $T^{\cn} M$ and
define the symmetric bilinear operator $\kappa$ by
$$\kappa(\phi,\psi)= g(\text{grad}(\phi),\text{grad}(\psi)).$$
Two maps $\phi,\psi: M\to\cn$ are said to be {\it orthogonal} if
$\kappa(\phi,\psi)=0$. The harmonicity and horizontal conformality
of $\phi:(M,g)\to\cn$ are then given by the following relations
$$\tau(\phi)=0\ \ \text{and}\ \ \kappa(\phi,\phi)=0.$$

\begin{definition}
Let $(M,g)$ be a semi-Riemannian manifold.  A set
$$\Omega=\{\phi_i:M\to\cn\ |\ i\in I\}$$ of complex-valued
functions is said to be an {\it orthogonal harmonic
family} on $M$ if for all $\phi,\psi\in\Omega$
$$\tau(\phi)=0\ \ \text{and}\ \ \kappa(\phi,\psi)=0.$$
\end{definition}

The following result shows that the elements of an orthogonal harmonic
family can be used to produce a variety of harmonic morphisms.

\begin{proposition}\cite{Gud-1}\label{prop:local-sol}
Let $(M,g)$ be a semi-Riemannian manifold and
$$\Omega=\{\phi_k:M\to\cn\ |\ k=1,\dots ,n\}$$ be a finite orthogonal
harmonic family on $(M,g)$.  Let $\Phi:M\to\cn^n$ be the map given
by $\Phi=(\phi_1,\dots,\phi_n)$ and $U$ be an open subset of
$\cn^n$ containing the image $\Phi(M)$ of $\Phi$. If
$$\tilde\F=\{F_i:U\to\cn\ |\ i\in I\}$$ is a family of holomorphic
functions then $$\F=\{\psi:M\to\cn\ |\ \psi=F(\phi_1,\dots ,\phi_n
),\ F\in\tilde\F\}$$ is an orthogonal harmonic family on $(M,g)$.
\end{proposition}

\section{Symmetric spaces}

Let $(G/K,g)$ be a Riemannian symmetric space of non-compact type,
where $G$ is a non-compact, semi-simple Lie group and $K$ a
maximal compact subgroup of $G$. Then the Killing form
$$B:\g\times\g\to\rn$$ of the Lie algebra $\g$ of $G$ induces a
bi-invariant semi-Riemannian metric $\hat g$ on $G$. Furthermore
it induces an orthogonal decomposition
$$\g=\k\oplus\p$$  of $\g$, where $\k$ is the Lie algebra of $K$.
The restriction of the Killing form to the orthogonal complement
$\p$ of $\k$ in $\g$ is positive definite and the natural
projection $\pi :G\to G/K$ is a Riemannian submersion with totally
geodesic fibres and hence a harmonic morphism by Theorem
\ref{theo:semi-B-E}.

Employing Proposition \ref{prop-lift} we see that the problem of
finding harmonic morphisms defined on an open subset $W$ of the
Riemannian symmetric space $G/K$ is equivalent to the problem of
finding $K$-invariant harmonic morphisms on the open subset $\hat
W=\pi^{-1}(W)$ of the Lie group $G$. If $Z$ is an element of the
Lie algebra $\g$ of left invariant vector fields on G and
$\phi:\hat W\to\cn$ is a map defined locally on $G$, then
$$Z(\phi)(p)=\frac {d}{ds}\big|_{s=0}\phi(p\cdot\exp(sZ)),$$
$$Z^2(\phi)(p)=\frac {d^2}{ds^2}\big|_{s=0}\phi(p\cdot\exp(sZ)).$$
If  in addition, the map $\phi$ is $K$-invariant and $Z\in\k$ then
$$Z(\phi)=0\ \ \text{and}\ \ Z^2(\phi)=0.$$
This means that if $\phi,\psi:\hat W\to\cn$ are complex-valued,
$K$-invariant maps defined locally on $G$, then the tension field
$\tau(\phi)$ and the $\kappa$-operator $\kappa(\phi,\psi)$ are
given by
$$\tau(\phi)=\sum_{Z\in\B}Z^2(\phi)\qquad\text{and}\qquad
\kappa(\phi,\psi)=\sum_{Z\in\B}Z(\phi)Z(\psi),$$ where $\B$ is any
orthonormal basis of the orthogonal complement $\p$ of $\k$ in
$\g$.

We now show how a locally defined complex-valued harmonic
morphism from a Riemannian symmetric space $G/K$ of non-compact
type gives rise to a locally defined harmonic morphism from the
compact dual space $U/K$ and vice versa. Recall that any harmonic
morphism between real analytic Riemannian manifolds is real
analytic, see \cite{Bai-Woo-book}.

Let $W$ be an open subset of $G/K$ and $\phi:W\to\cn$ be a real
analytic map. By composing $\phi$ with the natural projection
$\pi:G\to G/K$ we obtain a real analytic $K$-invariant map
$\hat\phi:\hat W\to\cn$ from the open subset $\hat W=\pi^{-1}(W)$
of $G$. Let $G^\cn$ denote the complexification of the Lie group
$G$. Then $\hat\phi$ extends uniquely to a $K$-invariant
holomorphic map $\phi^\cn:W^\cn\to\cn$ from some open subset
$W^\cn$ of $G^\cn$. By restricting this map to $U\cap W^\cn$ and
factoring through the natural projection $\pi^*:U\to U/K$ we
obtain a real analytic map $\phi^*:W^*\to\cn$ from some open
subset $W^*$ of $U/K$.

\begin{theorem}\label{theo:dual}\cite{Gud-Sve-1}
Let $\F$ be a family of real analytic maps $\phi:W\to\cn$ locally
defined on the non-compact irreducible Riemannian symmetric space
$G/K$ and $\F^*$ be the dual family consisting of the maps
$\phi^*:W^*\to\cn$ ´locally defined on the dual space $U/K$
constructed as above. Then $\F$ is an orthogonal harmonic family
on ¤W¤ if and only if $\F^*$ is an orthogonal harmonic family on
$W^*$.
\end{theorem}

\section{The symmetric space $\SLR n/\SO n$}

In this section we construct $\SO n$-invariant harmonic morphisms
on the special linear groups $\SLR n$ inducing globally defined
harmonic morphisms from the irreducible Riemannian symmetric
spaces $\SLR n/\SO n$.  This leads to a proof of Conjecture
\ref{conj:existence} in these cases.

Let $\GLRP n$ be the connected component of the general linear
group $\GLR n$ containing the identity element i.e the set of real
$n\times n$ matrices with positive determinant. On its Lie algebra
$\glr n$ we have a bi-linear form $$(X,Y)\mapsto\trace XY.$$ This
form induces a bi-invariant semi-Riemannian metric on $\GLRP n$.
We also get an orthogonal decomposition $$\glr n=\so n\oplus\p$$
of $\glr n$, where $$\so n=\{Y\in\glr n|\ Y+Y^t=0\},\quad
\p=\{X\in \glr n|\ X-X^t=0\}.$$ The restriction of the form to
$\p$ induces a $\GLRP n$-invariant metric on the quotient manifold
$\GLRP n/\SO n$ turning it into a Riemannian symmetric space. The
homogeneous projection
$$\GLRP n\to\GLRP n/\SO n$$ is a Riemannian submersion with
totally geodesic fibres, hence a submersive harmonic morphism. The
isomorphism
$$\rn^+\times\SLR n\to\GLRP n,\qquad (r,x)\mapsto r x,$$ induces
an isometry
$$\rn^+\times\SLR n/\SO n\cong\GLRP n/\SO n,$$ which is simply the
de Rham decomposition of $\GLRP n/\SO n$. Hence the problem of
finding harmonic morphisms on $\SLR n/\SO n$ is equivalent to
finding harmonic morphisms on $\GLRP n$ which are invariant under
the action of $\rn^+\times\SO n$.

\begin{theorem}\label{theo:real}
Let $\Phi,\Psi:\GLRP n\to\rn^{n\times n}$ be the $\SO n$-invariant
matrix valued maps, where
$$\Phi:x\mapsto x x^t,\quad \Phi=[\phi_{kl}]_{k,l=1}^n$$ and the
components $\psi_{kl}:\GLRP n\to\rn$ of $\Psi$ satisfy
$$\psi_{kl}=\sqrt{\phi_{kk}\phi_{ll}-\phi_{kl}^2}.$$
If $k\neq l$, then the map $$(\phi_{kl}+i\psi_{kl})/\phi_{ll}$$ is
a globally defined $\rn^+\times\SO n$-invariant harmonic morphism
on $\GLRP n$ inducing a globally defined harmonic morphism on the
irreducible Riemannian symmetric space $\SLR n/\SO n$.
\end{theorem}

For a proof of Theorem \ref{theo:real}, see Appendix
\ref{app:real}.

\section{The Symmetric Space $\SUS{2n}/\Sp n$}\label{sect:quat}

In this section we construct $\Sp n$-invariant harmonic morphisms
on the Lie groups $\SUS{2n}$ inducing globally defined harmonic
morphisms from the irreducible Riemannian symmetric spaces
$\SUS{2n}/\Sp n$. This leads to a proof of Conjecture
\ref{conj:existence} in these cases.

The quaternionic general linear group $\GLH n$ has a well-known
complex representation $$\US{2n}=\{z\in\GLC{2n}\ |\ zJ=J\bar z\},
$$ and its Lie algebra $\us{2n}$ is given by
$$\us{2n}=\{\begin{pmatrix}\alpha & \beta \\ -\bar\beta &
\bar\alpha\end{pmatrix} \ |\ \alpha, \beta\in\glc n\}.$$ The
quaternionic analogue to the complex special linear group $\SLC n$
is given by $$\SUS{2n}=\{z\in\SLC{2n}\ |\ zJ=J\bar z\},\text{
where }J =\begin{pmatrix}0 & I_n \\ -I_n & 0\end{pmatrix}$$ and
the compact Lie subgroup $\Sp n$ is defined by
$$\Sp{n}=\{z\in\SU{2n}\ |\ zJ=J\bar z\}=\SUS{2n}\cap\SU{2n}.$$
On the Lie algebra $\us{2n}$ we have the symmetric bi-linear form
$$(X,Y)\mapsto\Re(\trace XY)$$ inducing a bi-invariant
semi-Riemannian metric on the Lie group $\US{2n}$. This gives the
orthogonal decomposition
$$\us{2n}=\sp n\oplus\p$$ of $\us{2n}$, where $\sp n$ is the Lie algebra
\begin{eqnarray*}
\sp{n}&=&\{Y\in \us{2n}\ |\ Y^*+Y=0\}\\
&=&\{\begin{pmatrix}\alpha & \beta
\\ -\bar\beta & \bar\alpha\end{pmatrix}
\ |\ \alpha^*+\alpha=0,\ \beta^t-\beta=0\} \end{eqnarray*} of $\Sp
n$ and the orthogonal complement $\p$ is given by
\begin{eqnarray*}
\p&=&\{X\in\sus{2n}\ |\ X^*-X=0\}\\
&=&\{\begin{pmatrix} \alpha & \beta \\ -\bar\beta &
\bar\alpha\end{pmatrix} \ |\ \alpha^*-\alpha=0,\
\beta^t+\beta=0\}.
\end{eqnarray*}
The restriction of the form to $\p$ induces a $\US{2n}$-invariant
metric on the quotient manifold $\US{2n}/\Sp n$ turning it into a
Riemannian symetric space.  The homogeneous projection
$$\US{2n}\to\US{2n}/\Sp n$$ is a Riemannian submersion with
totally geodesic fibres, hence a submersive harmonic morphism.

The determinant of any element of $\US{2n}$ is a positive real
number. We thus get an isomorphism
$$\rn^+\times\SUS{2n}\to\US{2n},\qquad (r,z)\mapsto r z$$
inducing an isometry $$\rn^+\times\SUS{2n}/\Sp n\cong\US{2n}/\Sp
n,$$ which is simply the de Rham decomposition of $\US{2n}/\Sp n$.
Hence the problem of finding harmonic morphisms on the irreducible
Riemannian symmetric space $\SUS{2n}/\Sp n$ is equivalent to
finding harmonic morphisms on $\US{2n}$ which are invariant under
the action of $\rn^+\times\Sp n$.

\begin{theorem}\label{theo:quat}
Let $\Phi:\US{2n}\to\cn^{2n\times 2n}$ be the $\Sp n$-invariant
matrix-valued map given by $$\Phi:x\mapsto x x^*,\quad
\Phi=[\phi_{kl}]_{k,l=1}^{2n}.$$  For any $1\le l\le n$, the set
$$\F_l=\{\phi_{kl}/\phi_{ll}:\US{2n}\to\cn|\ k\neq l\}$$ is an
orthogonal harmonic family on $\US{2n}$ of $\rn^+\times\Sp
n$-invariant maps, inducing a globally defined orthogonal harmonic
family on the irreducible Riemannian symmetric space $\SUS{2n}/\Sp
n$.
\end{theorem}

For a proof of Theorem \ref{theo:quat}, see Appendix
\ref{app:quat}.

\section{The compact dual cases}

In this section we employ the duality principle of Theorem
\ref{theo:dual} to construct locally defined harmonic morphisms on
the compact irreducible Riemannian symmetric spaces
$\SU{n}/\SO{n}$ and $\SU{2n}/\Sp{n}$.  These are the dual spaces
to the non-compact $\SLR n/\SO n$ and $\SU{2n}/\Sp{n}$ studied in
the previous two sections.

\begin{theorem}
Let $\Phi^*,\Psi^*:\SU n\to\cn^{n\times n}$ be the $\SO
n$-invariant, matrix valued maps defined on the special unitary
group $\SU n$ where $$\Phi^*:x\mapsto x x^t,\quad
\Phi^*=[\phi^*_{kl}]_{k,l=1}^n$$ and the components
$\psi^*_{kl}:\SU n\to\cn$ of $\Psi^*$ satisfy
$$\psi^*_{kl}=\sqrt{\phi^*_{kk}\phi^*_{ll}-{\phi^*_{kl}}^2}.$$
Here the complex square root is the standard extension of the
classical real root. If $k\neq l$, then the $\SO n$-invariant,
complex-valued  map $$(\phi^*_{kl}+i\psi^*_{kl})/\phi^*_{ll}$$ is
an harmonic morphism defined on the open subset
$$\widehat W^*_{kl}=\{x\in\SU n\ |\ \phi^*_{ll}(x)\neq 0,\
\phi^*_{kk}(x)\phi^*_{ll}(x)-{\phi^*_{kl}(x)}^2\notin i\rn \}$$ of
$\SU n$, inducing a harmonic morphism locally defined on the open
subset $\widehat W^*_{kl}/\SO n$ of the irreducible Riemannian
symmetric space $\SU n/\SO n$.
\end{theorem}

\begin{proof} Let $\Phi$ and $\Psi$ be as defined in
Theorem \ref{theo:real}. Then the maps $\Phi^*$ and
$\Psi^*$ are restrictions to $\SU n$ of holomorphic maps
on $\GLC n$, the restrictions of which to $\GLRP n$
coincide with $\Phi$ and $\Psi$, respectively. The result
now follows from Theorem \ref{theo:real} and Theorem \ref{theo:dual}.
\end{proof}

For the spaces $\SU{2n}/\Sp n$ we have the following result.

\begin{theorem} Let $\Phi^*:\SU{2n}\to\cn^{2n\times 2n}$ be
the $\Sp n$-invariant matrix-valued map defined by
$$\Phi^*:x\mapsto x J^tx^tJ,\quad
\Phi^*=[\phi^*_{kl}]_{k,l=1}^{2n}.$$ For $1\le l\le n$, let $\S_l$
be the open subset $\{x\in\SU{2n}\ |\ \phi^*_{ll}\neq 0\}$ of
$\SU{2n}$. Then
$$\F^*_l=\{\phi^*_{kl}/\phi^*_{ll}:\S_k\to\cn\ |\ k\neq l\}$$
is an orthogonal harmonic family of $\Sp n$-invariant maps,
inducing a locally defined harmonic orthogonal family on the
irreducible Riemannian symmetric space $\SU{2n}/\Sp n$.
\end{theorem}

\begin{proof} Let $\Phi$ be as defined in Theorem \ref{theo:quat}.
Note that for $x\in\US{2n}$, $$\bar x=J^t x J.$$
Thus $\Phi^*$ is the restriction to $\SU{2n}$ of a holomorphic map
on $\GLC n$, the restriction of which to $\US{2n}$ coincides with $\Phi$.
The result now follows from Theorem \ref{theo:quat} and
Theorem \ref{theo:dual}.
\end{proof}

\section{symmetric spaces of type IV}

In this section we prove Conjecture \ref{conj:existence} for type IV Riemannian
symmetric spaces. These are of the form $G^\cn/G$, where $G$ is a compact simple
Lie group with complexification $G^\cn$. The dual space to $G^\cn/G$ is the type
II space $G$ with any bi-invariant metric. In \cite{Gud-Sve-2} we proved the
conjecture for the type II spaces, and so, by the duality principle, we know
that harmonic morphisms do exist \emph{locally} on any type IV symmetric space.
In the same paper we proved that these maps are \emph{globally} defined when
$G$ admits a Hermitian symmetric quotient, i.e. when $G$ is any compact simple
Lie group except $E_8$, $F_4$ or $G_2$.

We provide here a review of the construction of harmonic morphisms from type
II spaces and what the dual maps from type IV spaces look like. We also give
a self-contained proof that they can be chosen to be globally defined, with
no assumption on the existence of a Hermitian symmetric quotient.

For the sake of generality, let $G$ be a semisimple, compact Lie group with
complexification $G^\cn$ and let $\g$ and $\g^\cn$ be the corresponding Lie
algebras. Choose a maximal Abelian subalgebra $\h$ of $\g$ corresponding to
a maximal torus $H$ of $G$. The complexification $\h^\cn$ of $\h$ is then a
Cartan subalgebra of $\g^\cn$. Fix some ordering of the roots
$\Delta=\Delta^+\cup\Delta^-$ and denote by $\Pi$ the set of simple roots.
We obtain a Borel subalgebra
\begin{equation*}
\b=\h^\cn\oplus\sideset{}{^\oplus}\sum_{\alpha\in\Delta^+}\g_\alpha,
\end{equation*}
with corresponding Borel subgroup $B$ of $G^\cn$. It is well known that
the inclusion $G\hookrightarrow G^\cn$ induces a diffeomorphism
\begin{equation*}
G/H\cong G^\cn/B,
\end{equation*}
and either of these quotients is usually referred to as a \emph{(full)
flag manifold}, see e.g. \cite{Alek-Spiro}. It carries a complex structure and
the negative of the Killing form on $G$ induces a Hermitian,
cosymplectic metric, see \cite{Sve-2}. Thus, any local holomorphic
function on $G^\cn/B$ will be a harmonic morphism. Furthermore, any such
function $$\phi:U\subset G^\cn/B\cong G/H\to\cn$$ will lift to a locally
defined harmonic morphism on $G$. This shows that Conjecture
\ref{conj:existence} is true for any compact, semisimple Lie group,
in particular for any type II Riemannian symmetric space.

The corresponding \emph{dual} harmonic morphism, obtained by the
duality principle described in Theorem \ref{theo:dual}, is easily
seen (see \cite{Gud-Sve-2}) to be given by $$\phi^*:U^*\subset
G^\cn/G\to\cn,\qquad \phi^*(gG)=\phi(g\sigma(g)^{-1}B),$$ where
$\sigma$ is conjugation in $G^\cn$ with respect to $G$ and $$U^*
=\{gG\in G^\cn/G\ |\ g\sigma(g)^{-1}B\in U\}.$$ Thus, to prove
that $\phi^*$ can be chosen to be globally defined, we must show
that $U$ and $\phi$ can be chosen such that $$\phi:U\subset
G^\cn/B\to\cn$$ is non-constant and holomorphic, and $$\{gG\in
G^\cn/G\ |\ g\sigma(g)^{-1}B\in U\}=G^\cn/G.$$

Let $P_-$ be the nilpotent subgroup of $G^\cn$ with Lie algebra
\begin{equation*}
\p_-=\sum_{\alpha\in\Delta^+}\g_{-\alpha}.
\end{equation*}
The set $P_-B/B$ is often referred to as a
\emph{big cell}, i.e. open and dense in $G^\cn/B$, and biholomorphic
to $\cn^n$ for some $n$.

\begin{theorem}\label{theo:big-cell}
For any $g\in G^\cn$, the coset $g\sigma(g)^{-1}B\in G^\cn/B$
belongs to the big cell $P_-B/B$. Thus, for any non-constant
holomorphic function $$\phi:P_-B/B\to\cn,$$ the map
$$G^\cn/G\to\cn,\quad gG\mapsto\phi(g\sigma(g)^{-1}B)$$ is a
globally defined harmonic morphism on $G^\cn/G$. Moreover, such
non-constant holomorphic functions exist.
\end{theorem}

For each simple root $\alpha\in\Pi$, choose an element
$H_\alpha\in[\g_\alpha,\g_{-\alpha}]$ satisfying
$\alpha(H_\alpha)=2$. Associated to the ordering of the roots is
the (open) Weyl chamber
\begin{equation*}
\W=\{\lambda\in(\h^\cn)^*\ |\ \lambda(H_\alpha)>0\text{ for all }
\alpha\in\Pi\}.
\end{equation*}
By choosing an element $\lambda$ in the intersection of $\W$ and
the weight lattice, we obtain an irreducible representation $V$ of
$G^\cn$ with highest weight $\lambda$. On $V$ we fix a
$G$-invariant Hermitian product $\ip{\cdot}{\cdot}$; it follows
that the transpose $g^*$ of any element $g\in G^\cn$ equals
$\sigma(g)^{-1}$.

\begin{lemma}\label{lem:eigen} Any two weight spaces of $V$ are orthogonal.
\end{lemma}
\begin{proof} Assume that $\eta$ and $\gamma$ are two distinct
weights of $V$, and that $V_\eta$ and $V_\gamma$ are the corresponding
weight spaces. Take $v\in V_\eta$ and $w\in V_\gamma$. For any
$H\in\h$, $\eta(H)$ and $\gamma(H)$ are purely imaginary numbers. Thus
\begin{equation*}
\ip{v}{w}=\ip{\exp H\cdot v}{\exp H\cdot w}=\ip{e^{\eta(H)}v}{e^{\gamma(H)}w}
=e^{\eta(H)-\gamma(H)}\ip{v}{w}.
\end{equation*}
Hence we must have $\ip{v}{w}=0$.
\end{proof}

Denote by $\P V$ the projectivization of $V$. For any $w\in V\setminus\{0\}$,
denote by $[w]$ the corresponding element in $\P V$. Fix a non-zero vector
$v\in V_\lambda$ of highest weight. The group $G^\cn$ acts on $\P V$, and by
the particular choice of $\lambda$, the stabilizer of $[v]$ is precisely $B$.
Hence we have a realization of the flag manifold $G^\cn/B$ as the orbit
$G^\cn\cdot[v]$ in $\P V$.

Recall that the Weyl group $N_G(H)/H$ acts simply transitive on the set of Weyl
chambers, and also on the set of simple roots. Let $\omega$ be any representative 
of the unique element of the Weyl group taking $\Pi$ to $-\Pi$.

The main step in the proof of Theorem \ref{theo:big-cell} is the
following result, which is interesting in its own right.
\begin{proposition}\label{prop:big-cell} The image of the big cell
$P_-B/B$ in $\P V$ is equal to
$$\{[u]\in G^\cn\cdot[v]\ |\ \ip{u}{v}\neq 0\}.$$
\end{proposition}
\begin{proof} Since $$B\omega B\cdot[v]=B\omega\cdot[v]
=\omega\omega^{-1}B\omega\cdot[v]=\omega P_-\cdot[v],$$ it is
enough to prove that
$$B\omega\cdot[v]=\{[u]\in G^\cn\cdot[v]\ |\ \ip{u}{\omega\cdot
v}\neq0\}.$$

For any $b\in B$, note that
$$\ip{bw\cdot v}{w\cdot v}=\ip{w\cdot v}{b^*w\cdot v}
=\ip{w\cdot v}{ww^{-1}\sigma(b)^{-1}w\cdot v}.$$
Now $\sigma(b)^{-1}$ belongs to the Borel subgroup opposite to $B$
with Lie algebra $\p_-\oplus\h^\cn$, and $w^{-1}\sigma(b)^{-1}w\in B$.
Hence $w^{-1}\sigma(b)^{-1}w\cdot v$ is a non-zero multiple of $v$,
and so $$\ip{bw\cdot v}{w\cdot v}\neq0$$ for any $b\in B$. Next,
assume that $g\cdot v\in G^\cn\cdot[v]$ is such that
$$\ip{g\cdot v}{w\cdot v}\neq0.$$ According to the well-known Bruhat
decomposition of $G^\cn$, there is an element $\tilde w$ in the
Weyl group such that $g\in B\tilde w B$, i.e. we may write
$g=b\tilde w b'$ for some $b,b'\in B$. Then $b'\cdot v=xv$, for
some non-zero $x\in\cn$. Thus
$$0\neq\ip{g\cdot v}{w\cdot v}=\ip{b\tilde w b'\cdot v}{w\cdot v}
=x\ip{\tilde w\cdot v}{ww^{-1}\sigma(b)^{-1}w\cdot v}.$$
Since, as before, $w^{-1}\sigma(b)^{-1}w\cdot v$ is a non-zero
multiple of $v$, we conclude that $$\ip{\tilde w\cdot v}{w\cdot v}\neq0.$$
However, $\tilde w\cdot v\in V_{\tilde w(\lambda)}$ and
$w\cdot v\in V_{w(\lambda)}$. By Lemma \ref{lem:eigen}, we must have
$\tilde w(\lambda)=w(\lambda)$. As the Weyl group acts simply
transitively on the set of Weyl chambers we must have $\tilde w=w$
and thus $$g\cdot[v]=bw\cdot[v]\in Bw\cdot[v].$$
\end{proof}
\begin{proof}[Proof of Theorem \ref{theo:big-cell}] Since the big cell is
biholomorphic to $\cn^n$ for some $n$, it is clear that we can find a
non-constant holomorphic function $$\phi:P_-B/B\to\cn,$$ and this lifts
to a harmonic morphism, locally defined on $G$. In \cite{Gud-Sve-2} we
show that the dual map, locally defined on $G^\cn/G$, is given by
$$gG\to\phi(g\sigma(g)^{-1}B).$$ This is a harmonic morphism according
to Theorem \ref{theo:dual}. All that remains is to prove that this is
globally defined, i.e. that $$g\sigma(g)^{-1}B\in P_-B/B$$ for any
$g\in G^\cn$. This follows immediately from Proposition
\ref{prop:big-cell}, since
$$\ip{g\sigma(g)^{-1}\cdot v}{v}
=\ip{\sigma(g)^{-1}\cdot v}{\sigma(g)^{-1}\cdot v}\neq 0.$$
\end{proof}

\appendix
\section{}\label{app:real}

In this section we give a proof of Theorem \ref{theo:real}. For
this we introduce the following standard notation also employed in
Appendix \ref{app:quat}. For the positive integers $k,l$
satisfying $1\le k,l\le n$ we denote by $E_{kl}$ the elements of
the Lie algebra $\glr n$ given by
$$(E_{kl})_{ij}=\delta_{ki}\delta_{lj}$$ and by $D_k$ the diagonal
matrices $$D_k=E_{kk}.$$ For $1\le k<l\le n$ let $X_{kl}$ and
$Y_{kl}$ be the matrices satisfying
$$X_{kl}=\frac 1{\sqrt 2}(E_{kl}+E_{lk}), \qquad Y_{kl}=\frac
1{\sqrt 2}(E_{kl}-E_{lk}).$$  With this notation at hand we define
the orthonormal basis $$\B=\{D_k|\ 1\le k\le n\}\cup\{X_{kl}|\
1\le k<l\le n\}$$ for the subspace $\p$ of $\glr n$, induced by
the splitting $\glr n=\so n\oplus\p$.

\begin{lemma}\label{lemm:formula-real}
If $x,y,\alpha,\beta\in\cn^n$, then
$$\sum_{k<l}^n\alpha xX_{kl}y^tX_{kl}\beta^t
+\sum_{k=1}^n\alpha xD_ky^tD_k\beta^t=\frac{1}{2}(\alpha x^ty\beta^t
+y\alpha^tx\beta^t)$$ and
$$\sum_{k<l}^n\alpha xY_{kl}y^tY_{kl}\beta^t
=\frac{1}{2}(\alpha x^ty\beta^t-y\alpha^tx\beta^t).$$
\end{lemma}

\begin{proof}
Here we shall only prove the first equality and leave the second
for the reader as an exercise. Employing the notation introduced
above one easily shows that the following matrix equation holds
$$\sum_{k<l}^nxX_{kl}y^tX_{kl}
+\sum_{k=1}^nxD_ky^tD_k
=\bigg[\frac{x_iy_j+x_jy_j}{2}\bigg]_{i,j=1}^n.$$ This
obviously implies the following relation, proving the statement:
\begin{equation*}
\begin{split}
\sum_{k<l}^n\alpha xX_{kl}y^tX_{kl}\beta^t+\sum_{k=1}^n\alpha
xD_ky^tD_k\beta^t=&\frac{1}{2}\sum_{i,j=1}^n(\alpha_ix_iy_j\beta_j
+\alpha_ix_jy_i\beta_j)\\ =&\frac{1}{2}(\alpha
x^ty\beta^t+y\alpha^tx\beta^t).
\end{split}
\end{equation*}
\end{proof}

\begin{lemma}\label{lemm:real}
Let $\Phi,\Psi:\GLR n\to\rn^{n\times n}$ be the $\SO n$-invariant
matrix valued maps defined on the general linear group $\GLR n$
with $$\Phi:x\mapsto x x^t,\quad \Phi=[\phi_{kl}]_{k,l=1}^n$$ and
the components $\psi_{kl}:\GLR n\to\rn$ of $\Psi$ satisfy
$$\psi_{kl}=\sqrt{\phi_{kk}\phi_{ll}-\phi_{kl}^2}.$$ Then the
following relations hold:
\begin{itemize}
\item[(i)] $\tau(\phi_{kl})=2(n+1)\phi_{kl}$, \item[(ii)]
$\kappa(\phi_{kl},\phi_{ij})=2(\phi_{ki}\phi_{lj}+\phi_{kj}\phi_{li})$,
\item[(iii)] $\kappa(\phi_{kl},\psi_{kj})=2\phi_{kl}\psi_{kj}$,
\item[(iv)] $\kappa(\psi_{kl},\psi_{kl})=2\psi_{kl}^2$, \item[(v)]
$\tau(\psi_{kl})=2(n-1)\psi_{kl}$.
\end{itemize}
\end{lemma}

\begin{proof}
If $Z\in\p$ and $x\in\GLR n$, then differentiation of the
matrix-valued map $\Phi$ gives
$$Z(\Phi)=\frac{d}{ds}\big|_{s=0}x\exp(sZ)\exp(sZ^t)x^t
=\frac{d}{ds}\big|_{s=0}x\exp(2sZ)x^t=2xZx^t,$$ and
$$Z^2(\Phi)=\frac{d^2}{ds^2}\big|_{s=0}x\exp(2sZ)x^t=4xZ^2x^t.$$

(i) Summing over the basis $\B$ of the subspace $\p$ we
immediately obtain
$$\tau(\Phi)=4\sum_{Z\in\B}xZ^2x^t=4x(\sum_{Z\in\B}Z^2)x^t=2(n+1)\Phi.$$

(ii) Employing Lemma \ref{lemm:formula-real} and the above formula
for the first order derivatives of $\Phi$ we see that
\begin{eqnarray*}
\kappa(\phi_{kl},\phi_{ij})
&=&\sum_{Z\in\B}Z(\phi_{kl})Z(\phi_{ij})\\
&=&4\sum_{Z\in\B}\ip{x_kZ}{x_l}\ip{x_iZ}{x_j}\\
&=&x_i(\sum x_kZx_l^tZ)x_j\\
&=&2(\ip{x_k}{x_i}\ip{x_l}{x_j}+\ip{x_k}{x_j}\ip{x_l}{x_i}).
\end{eqnarray*}

(iii)-(iv) Differentiation of the identity
$\phi_{kl}^2+\psi_{kl}^2 =\phi_{kk}\phi_{ll}$ gives
$$2\psi_{kl}Z(\psi_{kl})
=\phi_{ll}Z(\phi_{kk})+\phi_{kk}Z(\phi_{ll})-2\phi_{kl}Z(\phi_{kl}).$$
The statements (iii)-(iv) are direct consequences of the
definitions of the operators $\tau,\kappa$, the result in (ii) and
the above formula for $Z(\psi_{kl})$.

(v) By differentiating the identity $\phi_{kl}^2+\psi_{kl}^2
=\phi_{kk}\phi_{ll}$ yet again we obtain
\begin{eqnarray*}
2\psi_{kl}Z^2(\psi_{kl})&=&-2Z(\psi_{kl})^2+Z^2(\phi_{ll})\phi_{kk}
+2Z(\phi_{ll})Z(\phi_{kk})\\
& &\quad +\phi_{ll}Z^2(\phi_{kk}) -2Z(\phi_{kl})-2\phi_{kl}Z^2(\phi_{kl}).
\end{eqnarray*}
Then using (i)-(iv) one easily obtains the statement of (v).
\end{proof}

\begin{proof}[Proof of Theorem \ref{theo:real}]
Let the functions $P,Q:\GLRP n\to\cn$ be defined by
$P=\phi_{kl}+i\psi_{kl}$ and $Q=\phi_{ll}$.  Employing Lemma
\ref{lemm:real} we see that
$$\kappa(P,Q)=2PQ+2\phi_{kl}Q,\ \  \kappa(P,P)=4P\phi_{kl},
\ \ \tau(P)=2(n-1)P+4\phi_{kl}.$$  Then the  basic relations
$$X(P/Q)=\frac{X(P)Q-X(Q)P}{Q^2},$$
$$X^2(P/Q)=\frac{Q^2X^2(P)-PQX^2(Q)-2QX(P)X(Q)+2PX(Q)X(Q)}{Q^3}$$
for the first and second order derivatives of the quotient $P/Q$
imply
$$Q^3\tau(P/Q)=[Q^2\tau(P)-PQ\tau(Q)-2Q\kappa(P,Q)+2P\kappa(Q,Q)],$$
$$Q^4\kappa(P/Q,P/Q)=[Q^2\kappa(P,P)-2PQ\kappa(P,Q)+P^2\kappa(Q,Q)].$$
The statement of Theorem \ref{theo:real} is then a direct
consequence of these last equations and Lemma \ref{lemm:real}.
\end{proof}

\section{}\label{app:quat}

In this section we give a proof Theorem \ref{theo:quat}. For this
we employ the notation $D_k, X_{kl}, Y_{kl}$ for the elements of
$\glr n$ introduced in Appendix \ref{app:real}.  As an orthonormal
basis $\B$ for the subspace $\p$, induced by the splitting
$\sus{2n}=\sp n\oplus\p$, we have the union of the following five
sets of matrices
\begin{equation*}
\begin{split}
\B_1&=\{\frac 1{\sqrt 2}\begin{pmatrix}D_{k} & 0 \\
0 & D_{k}\end{pmatrix}\ |\ 1\le k\le n\},\\
\B_2&=\{\frac 1{\sqrt 2}\begin{pmatrix}X_{kl} & 0 \\
0 & X_{kl}\end{pmatrix}\ |\ 1\le k<l\le n\},\\
\B_3&=\{\frac 1{\sqrt 2}\begin{pmatrix}iY_{kl} & 0 \\
0 & -iY_{kl}\end{pmatrix}\ |\ 1\le k<l\le n\},\\
\B_4&=\{\frac 1{\sqrt 2}\begin{pmatrix}0 & Y_{kl} \\
-Y_{kl} & 0\end{pmatrix}\ |\ 1\le k<l\le n\},\\
\B_5&=\{\frac 1{\sqrt 2}\begin{pmatrix}0 & iY_{kl} \\
iY_{kl} & 0\end{pmatrix}\ |\ 1\le k<l\le n\}.
\end{split}
\end{equation*}

\begin{lemma}\label{lem:long}  If $x,y,\alpha,\beta\in\cn^{2n}$, then
$$\sum_{Z\in\B}\ip{\alpha Z}{\beta}\ip{xZ}{y}
=\frac{1}{2}(\ip{x}{\beta}\overline{\ip{y}{\alpha}}
+\omega(x,\alpha)\overline{\omega(y,\beta)}),$$ where
$\ip{}{},\omega:\cn^{2n}\to\cn$ are the forms defined by
$$\ip{x}{y}=xy^*\ \ \text{ and }\ \ \omega(x,y)=xJy^t.$$
\end{lemma}

\begin{proof}
For $x\in\cn^{2n}$ let us write $x=(x_1,x_2)$ where
$x_1,x_2\in\cn^n$. Then we have
$$\sum_{Z\in\B}\ip{\alpha Z}{\beta}\ip{xZ}{y}
=x(\sum_{Z\in\B}\alpha Z\beta^*Z)y^*$$ and the right-hand side
satisfies
\begin{equation*}
\begin{split}
& x(\sum_{Z\in\B_1}\alpha Z\beta^*Z)y^*\\
=&\frac{1}{2}\sum_{k}x
\begin{pmatrix}\alpha_1D_k\beta_1^*D_k+\alpha_2D_k\beta_2^*D_k & 0 \\
0 & \alpha_1D_k\beta_1^*D_k+\alpha_2D_k\beta_2^*D_k\end{pmatrix}y^*\\
=&\frac{1}{2}\sum_{k}(x_1(\alpha_1D_k\beta_1^*D_k+\alpha_2D_k\beta_2^*D_k)
y_1^*\\
 &\qquad\qquad
+x_2(\alpha_1D_k\beta_1^*D_k+\alpha_2D_k\beta_2^*D_k)y_2^*),
\end{split}
\end{equation*}
\begin{equation*}
\begin{split}
& x(\sum_{Z\in\B_2}\alpha Z\beta^*Z)y^*\\
=&\frac{1}{2}\sum_{k<l}x
\begin{pmatrix}\alpha_1X_{kl}\beta_1^*X_{kl}
  +\alpha_2X_{kl}\beta_2^*X_{kl} & 0 \\
 0 & \alpha_1X_{kl}\beta_1^*X_{kl}
 +\alpha_2X_{kl}\beta_2^*X_{kl}\end{pmatrix}y^*\\
=&\frac{1}{2}\sum_{k<l}(x_1(\alpha_1X_{kl}\beta_1^*X_{kl}
+\alpha_2X_{kl}\beta_2^*X_{kl})y_1^*\\
&\qquad\qquad +x_2(\alpha_1X_{kl}\beta_1^*X_{kl}
+\alpha_2X_{kl}\beta_2^*X_{kl})y_2^*).
\end{split}
\end{equation*}
Adding the two relations and applying Lemma
\ref{lemm:formula-real} gives
\begin{equation*}
\begin{split}
&\sum_{Z\in\B_1\cup\B_2}\ip{\alpha Z}{\beta}\ip{xZ}{y}\\
=&\frac{1}{4}(x_1\alpha_1^t\overline{\beta_1y_1^t}
+x_1\beta_1^*\alpha_1y_1^*+x_2\alpha_1^t\overline{\beta_1y_2^t}
+x_2\beta_1^*\alpha_1y_2^*\\
&+x_1\alpha_2^t\overline{\beta_2y_1^t}
+x_1\beta_2^*\alpha_2y_1^*+x_2\alpha_2^t\overline{\beta_2y_2}
+x_2\beta_2^*\alpha_2y_2^*).
\end{split}
\end{equation*}
Similar calculations yield
\begin{equation*}
\begin{split}
x(&\sum_{Z\in\B_3}\alpha Z\beta^*Z)y^*+x(\sum_{Z\in\B_4}\alpha Z\beta^*Z)y^*
+x(\sum_{Z\in\B_5}\alpha Z\beta^*Z)y^*\\
=&\frac{1}{2}\sum_{k<l}x\begin{pmatrix}-\alpha_1Y_{kl}\beta^*Y_{kl}
+\alpha_2Y_{kl}\beta_2^*Y_{kl} & 0 \\ 0 & \alpha_1Y_{kl}\beta^*Y_{kl}
-\alpha_2Y_{kl}\beta_2^*Y_{kl}\end{pmatrix}y^*\\
&+\frac{1}{2}\sum_{k<l}x\begin{pmatrix}0 & \alpha_1Y_{kl}\beta_2^*Y_{kl}
-\alpha_2Y_{kl}\beta_1^*Y_{kl} \\ -\alpha_1Y_{kl}\beta_2^*Y_{kl}
+\alpha_2Y_{kl}\beta_1^*Y_{kl} & 0\end{pmatrix}y^*\\
&+\frac{1}{2}\sum_{k<l}x\begin{pmatrix}0 & -\alpha_1Y_{kl}\beta_2^*Y_{kl}
-\alpha_2Y_{kl}\beta_1^*Y_{kl} \\ -\alpha_1Y_{kl}\beta_2^*Y_{kl}
-\alpha_2Y_{kl}\beta_1^*Y_{kl} & 0\end{pmatrix}y^*\\
=&\frac{1}{4}(-x_1\alpha_1^t\overline{\beta_1y_1^t}+x_1\beta_1^*\alpha_1y_1^*
+x_1\alpha_2^t\overline{\beta_2y_1^t}-x_1\beta_2^*\alpha_2y_1^*
+x_2\alpha_1^t\overline{\beta_1y_2^t}\\
&-x_2\beta_1^*\alpha_1y_2^*-x_2\alpha_2^t\overline{\beta_2y_2^t}
+x_2\beta_2^*\alpha_2y_2^*)\\
&+\frac{1}{2}(-x_1\alpha_2^t\overline{\beta_1y_2^t}+x_1\beta_1^*\alpha_2y_2^*
-x_2\alpha_1^t\overline{\beta_2y_1^t}+x_2\beta_2^*\alpha_1y_1^*).
\end{split}
\end{equation*}
Adding up we finally get the following
\begin{equation*}
\begin{split}
\sum_{Z\in\B}\ip{\alpha Z}{\beta}\ip{xZ}{y}
=&\frac{1}{2}(x_1\beta_1^*\alpha_1y_1^*+x_2\alpha_1^t\overline{\beta_1y_2^t}
+x_1\alpha_2^t\overline{\beta_2y_1^t}+x_2\beta_2^*\alpha_2y_2^*\\
&-x_1\alpha_2^t\overline{\beta_1y_2^t}+x_1\beta_1^*\alpha_2y_2^*
-x_2\alpha_1^t\overline{\beta_2y_1^t}+x_2\beta_2^*\alpha_1y_1^*)\\
=&\frac{1}{2}(\ip{x}{\beta}\overline{\ip{y}{\alpha}}
+\omega(x,\alpha)\overline{\omega(y,\beta)}).
\end{split}
\end{equation*}
\end{proof}

\begin{lemma}
Let $\Phi:\US{2n}\to\cn^{2n\times 2n}$ be the $\Sp n$-invariant
matrix-valued map defined by $$\Phi:x\mapsto x x^*,\quad
\Phi=[\phi_{kl}]_{k,l=1}^{2n}.$$ Then the following relations hold
\begin{itemize}
\item[(i)] $\tau(\phi_{kl})=(4n-2)\phi_{kl}$, \item[(ii)]
$\kappa(\phi_{kl},\phi_{rl})=2\phi_{kl}\phi_{rl}$.
\end{itemize}
\end{lemma}

\begin{proof}
If $Z\in\p$ and $x\in\US{2n}$, then differentiation of the matrix
valued map $\Phi$ gives
$$Z(\Phi)=\frac{d}{dt}\big|_{t=0}x\exp(tZ)\exp(tZ^*)x^*
=\frac{d}{dt}\big|_{t=0}x\exp(2tZ)x^*=2xZx^*,$$ and
$$Z^2(\Phi)=\frac{d^2}{dt^2}\big|_{t=0}x\exp(2tZ)x^*=4xZ^2x.$$
By summing over the basis $\B$ of the subspace $\p$ we immediately
get
$$\tau(\Phi)=4\sum_{Z\in\B}xZ^2x^*=4x(\sum_{Z\in\B}Z^2)x^*=(4n-2)\Phi.$$
A simple calculation and Lemma \ref{lem:long} now show that
\begin{equation*}
\kappa(\phi_{kl},\phi_{ij})=4\sum_{Z\in\B}\ip{x_kZ}{x_l}\ip{x_iZ}{x_j}
=2(\phi_{il}\phi_{kj}+\omega(x_i,x_k)\overline{\omega(x_j,x_l)}).
\end{equation*}
The second statement is an immediate consequence of the fact that
the form $\omega$ is skew-symmetric.
\end{proof}

\begin{proof}[Proof of Theorem \ref{theo:quat}]
Employing the chain rule we ´see that
$$\tau(\phi_{kl}/\phi_{ll})=\frac{\phi_{ll}^2\tau(\phi_{kl})
-\phi_{kl}\phi_{ll}\tau(\phi_{ll})-2\phi_{ll}\kappa(\phi_{kl},\phi_{ll})
+2\phi_{kl}\kappa(\phi_{ll},\phi_{ll})}{\phi_{ll}^3},$$ and
\begin{eqnarray*}
&&\kappa(\phi_{kl}/\phi_{ll},\phi_{rl}/\phi_{ll})\\
&=&\frac{\phi_{ll}^2\kappa(\phi_{kl},\phi_{rl})
-\phi_{rl}\phi_{ll}\kappa(\phi_{kl},\phi_{ll})
-\phi_{kl}\phi_{ll}\kappa(\phi_{rl},\phi_{ll})
+\phi_{kl}\phi_{rl}\kappa(\phi_{ll},\phi_{ll})}{\phi_{ll}^2}.
\end{eqnarray*}
Combining the results of Lemma \ref{lem:long} and the above
equations, we easily obtain the statement of Theorem
\ref{theo:quat}.
\end{proof}

\end{document}